\documentclass{amsart}

\usepackage{amsfonts}
\usepackage{amsmath}
\usepackage{amsthm}
\usepackage{graphicx}

\newtheorem{pro}{Proposition}[section]
\newtheorem{thm}[pro]{Theorem}
\newtheorem{lem}[pro]{Lemma}

\newtheorem{example}[pro]{Example}

\newtheorem{rmkks}[pro]{Remarks}

\newtheorem{cor}[pro]{Corollary}

\theoremstyle{definition}
\newtheorem{dfn}[pro]{Definition}
\newtheorem{dfns}[pro]{Definitions}

\newtheorem*{note}{Note}

\newcommand{\del}{\partial}

\newcommand{\ie}{{\it i.e.}}

\newcommand{\ess}{essential}

\newcommand{\hhs}{Heegaard surface}
\newcommand{\hhsp}{Heegaard splitting}

\newcommand{\nbhd}{neighborhood}

\newcommand{\scc}{simple closed curve}
\newcommand{\sccs}{simple closed curves}

\newcommand{\st}{solid torus}

\newcommand{\mx}{M_X}
\newcommand{\mmx}{M \setminus X}

\newcommand{\s}{\Sigma}
\newcommand{\smx}{\s_{M \setminus X}}
\newcommand{\sx}{\s_X}

\title[Local detection of strongly irreducibility]
{Local detection of strongly irreducible Heegaard splittings via
knot exteriors}

\date{\today}

\address{Department of Mathematics, Nara Women's University
Kitauoya Nishimachi, Nara 630-8506, Japan}
\address{Department of Mathematics, Nara Women's University
Kitauoya Nishimachi, Nara 630-8506, Japan and
Department of mathematical Sciences, University of
Arkansas, Fayetteville, AR 72701}

\email{tsuyoshi@cc.nara-wu.ac.jp}
\email{yoav@uark.edu}

\author{Tsuyoshi Kobayashi}
\author{Yo'av Rieck}

\thanks{The second author was supported by grant JSPS 00024}

\begin{document}
\maketitle

\section{Introduction}
\label{sec:intro}

Let $T$ be a compressible torus in an irreducible 3-manifold $M$
other than $S^3$. It is easy to see that either : 1. $T$ bounds a
solid torus, or: 2. $T$ bounds a submanifold homeomorphic to the
exterior of a non-trivial knot in $S^3$, where the compressing
disk for $T$ is a meridian disk of the knot.

The intersection of a strongly irreducible Heegaard surface with
solid tori was analyzed
by Y.Moriah and H.Rubinstein \cite{mr} and M.Scharlemann \cite{MR99h:57040},
where it was shown that such intersection can only occur in a very
restricted manner. The purpose of this paper is to analyze the
intersection of a strongly irreducible Heegaard surface with the
knot exterior in the situation 2 above. Precisely, let $M$ be a
3-manifold with a strongly irreducible Heegaard splitting $H_1
\cup_\s H_2$. Let $X$ be a 3-dimensional submanifold such that:

\begin{enumerate}

\item $X$ is homeomorphic to the exterior of a non-trivial knot in
$S^3$, and---

\item there is a compressing disk, say $D_X$, of $\partial X$ such
that $\partial D_X$ is a meridian curve of $X$.

\end{enumerate}

\begin{figure}[ht]
\begin{center}
\includegraphics[width=7cm, clip]{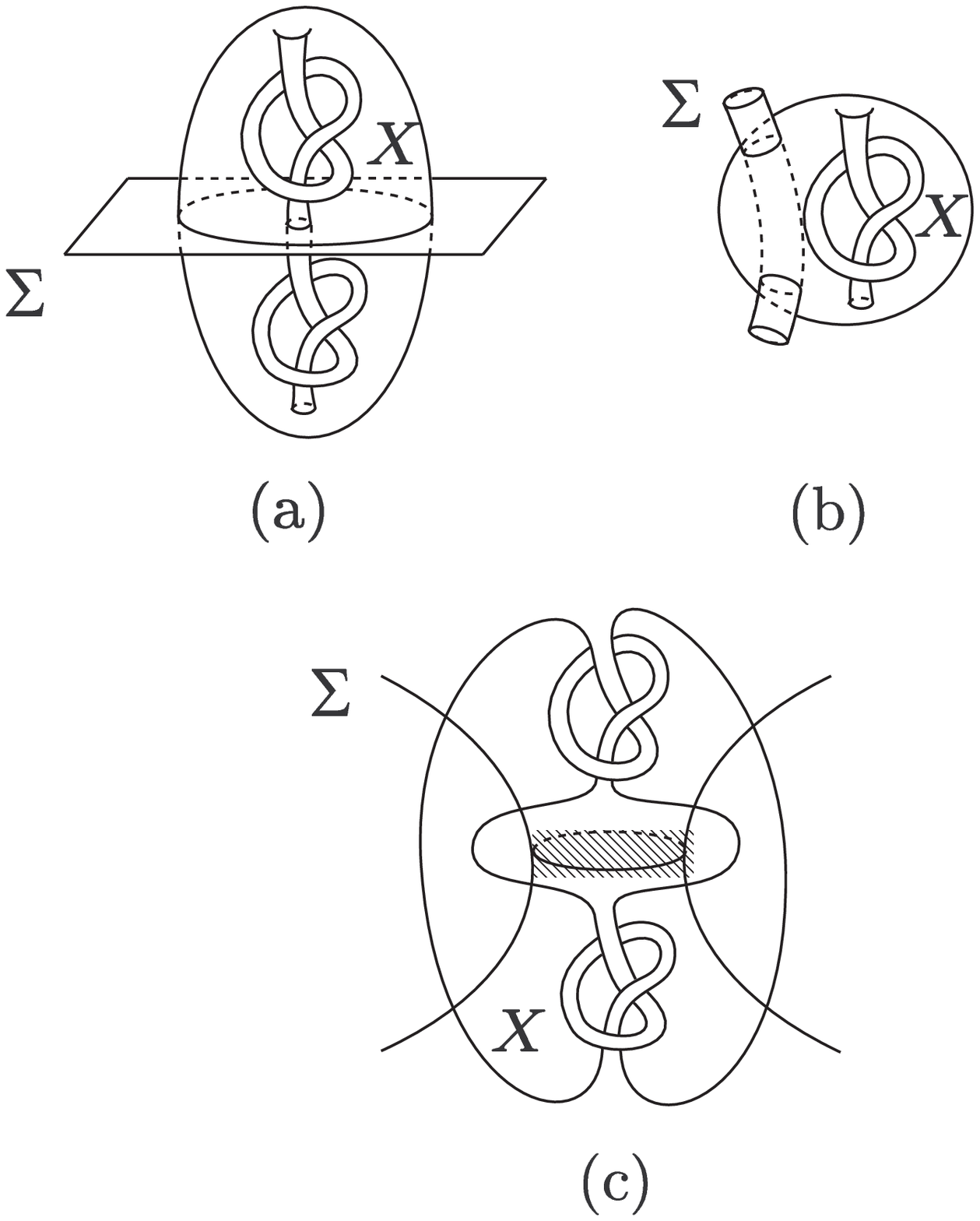}
\end{center}

\caption{}
\label{fig:abc}
\end{figure}

Note that $N(X \cup D_X)$ is a 3-ball, hence $X$ embeds in any
manifold, and there are several ways $X$ can intersect $\s$. These
are shown in Figure \ref{fig:abc}. In Figure \ref{fig:abc}(a) $\s$
intersects $T$ in \sccs\ which are essential in $T$ but
inessential in $\s$, and in Figure \ref{fig:abc}(b) $\s$
intersects $T$ in \sccs\ which are essential in $\s$ but
inessential in $T$. More interesting is Figure \ref{fig:abc}(c),
where all curves of intersection are \ess\ in both $T$ and $\s$.
(The part of the \hhs\ shown there is a cylinder, which is a
\nbhd\ of the boundary of the shaded meridian disk of $H_1$.) Note
that in Figure \ref{fig:abc}(c), the slope of $\s \cap T$ is
meridional, and each component of $\s \cap X$ is an annulus. We
call such an annulus a {\it meridional annulus}. A meridional
annulus in $X$ is either boundary parallel, or a decomposing
annulus in the exterior of a composite knot.

The main result of this paper is as follows.

\begin{thm}
\label{main theorem} Let $M$ be a 3-manifold other than the
3-sphere $S^3$ with a strongly irreducible Heegaard splitting $H_1
\cup_\s H_2$. Let $X$ be a 3-dimensional submanifold of $M$ such
that:

\begin{enumerate}

\item $X$ is homeomorphic to the exterior of a non-trivial knot in
$S^3$, and---

\item
there is a compressing disk, say $D_X$, of $\partial X$
such that $\partial D_X$ is a meridian curve of $X$.

\end{enumerate}

Suppose that $\partial X \cap \s$ consists of a non-empty
collection of \sccs\ which are essential in both $\partial X$ and
$\s$. Then we have:

\begin{enumerate}

\item the closure of some component of $\s \setminus
\partial X$ is an annulus and is parallel to an annulus in $\partial
X$, and---

\item each component of $\s \cap X$ is a (possibly boundary
parallel) meridional annulus.

\end{enumerate}

\end{thm}

\begin{rmkks}

\begin{enumerate}

\item The annulus in conclusion~1 of Theorem~\ref{main theorem}
may be contained in $X$ or in $\mbox{cl}(M \setminus X)$.

\item In \cite{mr} and \cite{MR99h:57040}, the
intersection of a strongly irreducible \hhs\ with a ball was considered,
and it is shown that if the boundary of the ball is incompressible in the
handlebodies then the intersection of the \hhs\ with the ball is
an unknotted planar surface.

\end{enumerate}
\end{rmkks}

Before finishing this section, we bring its main application. It
is concerned with a \em pair \em of strongly irreducible \hhs s
that intersect essentially and spinally.  In \cite{rs1:1996.1}
J.H.Rubinstein and Scharlemann showed that such intersection can
always be obtained, if we allow a single trivial \scc.  In
\cite{rieck:rubinstein} Rieck and Rubinstein show that the trivial
\scc\ can be avoided.

\begin{cor}
\label{cor:bounding-solid-torus}
Let $M$ be an irreducible, a-toroidal manifold, let
$\s_1,\;\s_2 \subset M$ be strongly irreducible \hhs s that
intersect essentially and spinally.

If $T \subset \s_1 \cup \s_2$ is a torus then $T$ bounds a \st.
\end{cor}

We prove the corollary assuming the theorem:

\begin{proof}[Proof of Corollary \ref{cor:bounding-solid-torus}]
Suppose $T$ does not bound a solid torus, then $T$ bounds a knot
exterior $X$ as in Theorem \ref{main theorem}.  Since $\s_1 \cup
\s_2$ is a finite complex we may pass to an innermost
counterexample to the corollary, \ie\ we may assume there does not
exist $T' \subset \s_1 \cup \s_2$ bounding a non-trivial knot
exterior $X'$ so that
$X'$ is strictly contained in $X$.
We shrink $X$ slightly to obtain the knot exterior $\widehat X$
and the torus $\widehat T = \del \widehat X$, so that $\widehat T$
is transverse to $\s_1$ and $\s_2$.  By essentialilty all the
curves of $\s_1 \cap \s_2$ are essential in $\s_1$, and since any
curve of $\s_1 \cap \widehat T$ is parallel to some such curve in
$\s_1$, it must also be essential in $\s_1$.
Furthermore, if a curve of $\s_1 \cap \widehat T$ is inessential
in $\widehat T$ it is parallel to a curve of $\s_1 \cap \s_2$ on
$T$ that is inessetial there, contradicting essentiallity. Hence
the conditions of Theorem \ref{main theorem} are satisfied. If
$\s_1 \cap \text{int}(X) \neq \emptyset$ any component of that
intersection yields (by Theorem \ref{main theorem}) a meridional
annulus in $\widehat X$. Since an annulus that decomposes a
non-trivial knot exterior into two solid tori is not meridional we
can use the meridional annulus and an annulus from $T$ to get a
torus $T'$ that contradicts our choice of $X$. ($T'$ would not
bound a solid torus on either side: on the side contained in $B$
as we just saw, and on the other side it bounds a piece in which a
punctured copy of $M$ is embedded.) Hence $\s_1 \cap \text{int}(X)
= \emptyset$ and similarly $\s_2 \cap \text{int}(X) = \emptyset$.
Thus $X$ is a component of $M$ cut open along $\s_1 \cup \s_2$ but
in \cite{rieck-proc} Rieck showed that every such component is a
handlebody (it is here that we use the spinality assumption), a
contradiction.
\end{proof}

\section{Preliminaries}
\label{sec:preliminaries}

Throughout this paper, we work in the differentiable category. For
a submanifold $H$ of a manifold $M$, $N(H,M)$ denotes a regular
neighborhood of $H$ in $M$. When $M$ is well understood we
abbreviate $N(H,M)$ to $N(H)$.  For the definitions of standard
terms in 3-dimensional topology, we refer to \cite{hempel} or
\cite{jaco}.

A 3-manifold $C$ is a {\it compression body} if there exists a
compact, connected (not necessarily closed) surface $F$ such that
$C$ is obtained from $F \times [0,1]$ by attaching 2-handles along
mutually disjoint simple closed curves in $F \times \{ 1 \}$ and
capping off the resulting 2-sphere boundary components which are
disjoint from $F \times \{ 0 \}$ by 3-handles. The subsurface of
$\partial C$ corresponding to $F \times \{ 0 \}$ is denoted by
$\partial_+C$. Then $\partial_-C$ denotes the subsurface
$\text{cl} (\partial C -(\partial F \times [0,1] \cup
\partial_+C))$ of $\partial C$. A compression body $C$ is called a
{\it handlebody} if $\partial_-C = \emptyset$. A compressing disk
$D (\subset C)$ of $\partial_+ C$ is called a {\it meridian disk}
of the compression body $C$.

\begin{rmkks}\label{remark of handlebody}
The following properties are known for
handlebodies

\begin{enumerate}

\item
Let $F$ be an incompressible surface in a handlebody.
Then either $F$ is boundary compressible or a meridian disk.

\item Let $F$ be an incompressible surface in a solid torus (\ie\
genus one handlebody). Then $F$ is either a meridian disk or a
boundary parallel annulus.

\item
Every incompressible surface in a handlebody cuts the handlebody
into handlebodies.

\end{enumerate}
\end{rmkks}

Let $N$ be a cobordism rel~$\partial$ between two surfaces $F_1$,
$F_2$ (possibly $F_1 = \emptyset$ or $F_2 = \emptyset$), i.e.,
$F_1$ and $F_2$ are mutually disjoint surfaces in $\partial N$
with $\partial F_1 \cong \partial F_2$ such that $\partial N = F_1
\cup F_2 \cup (\partial F_1 \times [0, 1])$, and $F_i \cap
(\partial F_1 \times [0, 1]) = \del F_i$ $(i=1, 2)$.

\begin{dfn}\label{Heegaard splitting}
We say that $C_1 \cup_P C_2$ (or $C_1 \cup C_2)$
is a {\it Heegaard splitting}
of $(N, F_1, F_2)$ (or simply, $N$) if it satisfies the following
conditions.

\begin{enumerate}
\item
$C_i$ $(i=1,2)$ is a compression body in $N$ such that
$\partial_- C_i = F_i$,

\item
$C_1 \cup C_2 = N$, and

\item
$C_1 \cap C_2 = \partial_+ C_1 = \partial_+ C_2 = P$.

\end{enumerate}

The surface $P$ is called a {\it Heegaard surface} of
$(N, F_1, F_2)$ (or, $N$).
\end{dfn}

\begin{dfns}\label{reducible Heegaard splitting}

\
\begin{enumerate}

\item
A Heegaard splitting $C_1 \cup_P C_2$ is {\it reducible}
if there exist meridian disks $D_1$, $D_2$ of
the compression bodies $C_1$, $C_2$
respectively such that
$\partial D_1= \partial D_2$

\item
A Heegaard splitting $C_1 \cup_P C_2$ is {\it weakly reducible}
if there exist meridian disks $D_1$, $D_2$ of
the compression bodies $C_1$, $C_2$
respectively such that
$\partial D_1 \cap  \partial D_2 = \emptyset$.
If $C_1 \cup_P C_2$ is not weakly reducible, then it is called
{\it strongly irreducible}.
\end{enumerate}
\end{dfns}

A {\it spine} of a handlebody $H$ is a 1-complex $L$ embedded in
$\text{int}H$ such that $L$ is a deformation retract of $H$. A
{\it cycle} of the spine $L$ is a \scc\  embedded in $L$. Then the
following is proved by C.Frohman \cite{frohman}, and will be used
in the proof of Theorem~\ref{main theorem}.

\begin{lem}[Frohman's Lemma]\label{Frohman}
Let $H_1 \cup_S H_2$ be a \hhsp\  of a closed irreducible 3-manifold
$M$, and $Y$ a spine of $H_j$ $(j=1$ or $2)$. Suppose that there
is a 3-ball $B^3$ in $M$ such that some cycle of $Y$ is contained
in $B^3$. Then $H_1 \cup_S H_2$ is reducible.
\end{lem}

The next lemma proved by Scharlemann \cite{MR99h:57040} is also
used in the proof.

\begin{lem}[No Nesting Lemma]\label{no nesting}
Suppose that $H_1 \cup_S H_2$ is a strongly irreducible \hhsp\ of
a 3-manifold $M$, and $F$ a disk in $M$ transverse to $S$ with
$\partial F \subset S$. Then $\partial F$ also bounds a disk in
$H_j$ $(j = 1$ or $2)$.
\end{lem}

\begin{example}
{\rm The $(1,1)$ curve on the standard torus in $S^3$ shows that
the transversality assumption is needed.}
\end{example}

The next lemma must be well known, but for the convenience of the
reader, we bring it with a proof.

\begin{lem}\label{boudary compressible implies compressible}
Let $N$ be a 3-manifold with a toral boundary component $T$. Let
$S$ be a 2-sided surface properly embedded in $N$ such that $S \cap T$
consists of \ess\ \sccs\ in $T$. Suppose that there is a boundary
compressing disk $\Delta$ for $S$ such that $\Delta$ compresses
$S$ into $T$, \ie , $\Delta \cap \partial N = \partial \Delta \cap
T$ is an arc, say $\alpha$, and $\Delta \cap S = \partial \Delta
\cap S$ is an \ess\ arc in $S$, say $\beta$, such that $\alpha \cup
\beta = \partial \Delta$. Then we have have either one of the following.

\begin{enumerate}

\item
$S$ is compressible. Moreover if $S$ is separating in $N$, then the
compression occurs in the same side as the boundary compression.

\item $S$ is an annulus; moreover, when $N$ is irreducible, $S$ is
boundary parallel.

\end{enumerate}

\end{lem}

\begin{proof}

\medskip
\noindent Case 1. $\partial \alpha$ is contained in a single
component, say $\ell$, of $S \cap T$.

\medskip
Since $S$ is 2-sided, neighborhoods of both endpoints of $\alpha$
are contained in the same side of $\ell$. Then there is a subarc,
say $\alpha'$, of $\ell$ such that $\alpha \cup \alpha'$ bounds a
disk $D$ in $T$. We isotope $\Delta$ by moving $\alpha$ to
$\alpha'$ along $D$ to obtain $\Delta'$ such that $\partial
\Delta' \subset S$. Since $\Delta \cap S$ is an \ess\ arc in $S$,
we see that $\Delta'$ is a compressing disk for $S$, and this
gives the conclusion~1.

\medskip
\noindent
Case 2. $\partial \alpha$ is contained in different components, say
$\ell_1$ and $\ell_2$, of $S \cap T$.

\medskip
Let $A (\subset T)$ be the annulus bounded by $\ell_1 \cup \ell_2$
such that $\alpha \subset A$. Let $D$ be a disk obtained from $A$
by boundary compressing along $\Delta$, hence $\partial D \subset
S$. If $\partial D$ is essential in $S$, then we have the
conclusion~1. If $\partial D$ bounds a disk in $S$, then we see
that $S$ is an annulus.  If in addition $N$ is irreducible, the
sphere obtained by compressing $S \cup A$ along $\Delta$ bounds
a ball, and we easily see that $S$ is boundary parallel;
hence conclusion~2 holds.

\end{proof}

\begin{dfn}
\label{dfn:essential in HB} A surface properly embedded in a
handlebody is called \em essential \em if it is incompressible and
not boundary parallel.
\end{dfn}


%
%
%
%

\begin{lem}\label{essential annulus preserves spine}
Let $A$ be separating essential annulus properly embedded in
a handlebody $H$.
Then there is a spine $Y$ of $H$ such that
$Y$ intersect $A$ in one point, and that
each component of $Y \setminus A$ contains a cycle of $Y$.
\end{lem}

\begin{proof}
By 1 of Remarks~\ref{remark of handlebody} $A$ is boundary
compressible, and let $D$ be a disk obtained from $A$ by a
boundary compression. Since $A$ is essential, $D$ is a meridian
disk of $H$. Since $A$ is separating in $H$, $D$ is also
separating in $H$. Hence we can find a spine $Y$ of $H$ such that
$Y$ intersects $D$ in one point, and that each component of $Y
\setminus D$ contains a cycle of $Y$. Note that $A$ is recovered
from $D$ by adding a band. We may suppose that the band is
disjoint from $Y$, hence $Y$ gives the conclusion of the lemma.
\end{proof}

\section{Proof of Theorem~\ref{main theorem}}
\label{sec:proof}

Let $M$, $H_1 \cup_\s H_2$, $X$, and $D_X$ be as in
Theorem~\ref{main theorem}. Let $T = \partial X$, $B=N(X \cup
D_X)$, and $\mx = \text{cl}( \mmx )$.

\begin{figure}[ht]

\begin{center}
\includegraphics[width=6cm, clip]{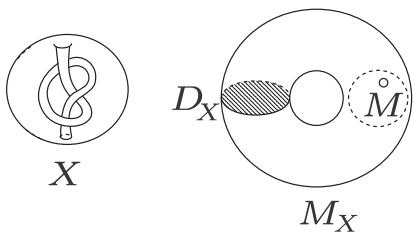}
\end{center}

\caption{$X$ and $\mx$.}
\label{fig:decomposing-m}
\end{figure}

Note that $B$ is a 3-ball in $M$, which contains $X$.
Note also that $\mx \cong (D^2 \times S^1) \# M$,
where the sphere $\partial B$ defines the connect sum structure,
and the disk $D_X$ is a meridian disk for $D^2 \times S^1$.
See Figure~\ref{fig:decomposing-m}.
We always assume $\del B \cap D_X = \emptyset$.
Recall that $X (\subset B)$ is in fact a
knot exterior in $S^3$ and the slope defined by $\del D_X$ on its
boundary is the slope of the trivial filling. We refer to $X$ as
\em the knot exterior. \em The slope of $\del D_X$ plays a crucial
role in our game and is called the {\it meridian} slope; $D_X$ is
called the meridian disk. Any other slope on $T$ is called {\it
longitudinal} if it intersects the meridional slope once, {\it
cabled} otherwise. Finally, we note that since $X$ is (by
assumption) a non-trivial knot exterior, $\del X$ is
incompressible in $X$, and on the boundary of $(D^2 \times S^1) \#
M$ only one slope compresses. Thus $D_X \subset M$ is the unique
compressing disk for $T$ (up-to isotopy relative to $T$), and the
only slope that compresses is the meridional slope.

\begin{proof}[Proof of Theorem \ref{main theorem}]

We divide the proof of the theorem into three steps. The first
(and main) step is:

\medskip
\noindent
{\bf Step 1: the slope of $\s \cap T$ is meridional.}

\medskip
Assume, for contradiction, that the slope is not meridional.
Note that each component of $T \cap H_i$ $(i=1,2)$ is an annulus.

\medskip
\noindent
\underline{Claim 1.} Each component of $T \cap H_i$
$(i=1,2)$ is incompressible in $H_i$.

\begin{proof}
Assume that there is a component $A$ of $T \cap H_j$ $(j=1$ or $2)$ such that
$A$ is compressible in $H_j$.
By using innermost disk arguments, we may suppose that
$\text{int} D \cap T = \emptyset$.
This shows that $D \subset \mx$, and $\partial D$ is a meridional slope.
Hence the slope of $\s \cap T$ is meridional,
contradicting the assumption of the proof of Step~1.
\end{proof}

\noindent
\underline{Claim 2.} By applying an isotopy, if
necessary, we may suppose that no component of $\s$ cut along $T$
is an annulus which is boundary parallel in $X$ or $\mx$.

\begin{proof}
Suppose there is such a component. Using it to guide an isotopy of
$\s$ we reduce $|\s \cap T|$ by two. Repeat the procedure as much
as possible. If we come to the situation that $\s \cap T \ne
\emptyset$, and no component of $\s$ cut along $T$ is boundary
parallel, then we are done. Assume that the intersection $\s \cap
T$ becomes empty. Then $\s$ is pushed into $X$ or $\mx$. However
the former is absurd ($\s$ is contained in the 3-ball $B$). Hence
$\s$ is pushed into $\mx$. Note that prior to the last isotopy $\s
\cap T$ consists of two simple closed curves, and we analyze this
configuration. Then $T \cap H_i$ $(i=1,2)$ consists of an annulus,
say $A_i$, and $A_j$ $(j=1$ or $2)$ is boundary parallel in $H_j$.
Since the argument is symmetric, we may suppose that $A_1$ is
boundary parallel in $H_1$.

\medskip
\noindent
Subclaim 2.1. $A_2$ is not boundary parallel in $H_2$.

\begin{proof}
Assume that $A_2$ is boundary parallel in $H_2$. Then either $T$
bounds a solid torus (if $A_1$ and $A_2$ are parallel to the same
annulus in $\s$), or $T$ is isotopic to $\s$ (if $A_1$ and $A_2$
are parallel to different annuli in $\s$), contradiction either
way.
\end{proof}

This together with Claim~1 shows that
$A_2$ is an essential annulus in $H_2$ and by
Lemma~\ref{essential annulus preserves spine} there is a cycle of a
spine of $H_2$ on each side of it.  But $A_2$ seperates $H_2$ into $X \cap
H_2$ and $\mx \cap H_2$ and so one of these cycles is contained in
$X$ and hence in $B$, and by Frohman's Lemma (\ref{Frohman}) $\s$
reduces, contradiction.


This completes the proof of Claim~2.
\end{proof}

\begin{note}
The argument in the proof of Claim~2 is a warm-up case of the
proof of Step~1, where we drive to find a cycle in $X$
(and hence in $B$) violating Frohman's Lemma.
\end{note}

\noindent
Notation: We denote $\s \cap \mx$ by $\smx$ and $\s \cap
X$ by $\sx$.

\medskip
\noindent
\underline{Claim~3.} By retaking $X$, if necessary, we
may suppose that no component of $\smx$ is an annulus.

\begin{proof}
Suppose that there is an annulus component, say $A$, in $\smx$.

\medskip
\noindent
Subclaim 3.1. $A$ is incompressible in $\mx$.

\begin{proof}
Assume that $A$ compressible in $\mx$. By compressing $A$, we
obtain two compressing disks for $\partial \mx$. This shows that
$\partial A$ is a meridional slope of $X$, contradicting the
assumption of the proof of Step~1.
\end{proof}

Recall that $\mx \cong (D^2 \times S^1) \# M$. Hence, by
Subclaim~3.1, 2 of Remarks~\ref{remark of handlebody}, and
Claim~2, we see that $A$ together with an annulus in $T$, say
$A'$, bounds a piece $P$ homeomorphic to $(D^2 \times S^1) \# M$,
where the slope of $\partial A$ is longitudinal in that solid
torus. Consider a torus, say $T'$, obtained by slightly pushing
$\partial P$ $(= A \cup A')$ into $P$. Let $P'$ be the submanifold
bounded by $T'$ which is contained in $P$.

\medskip
\noindent
Subclaim 3.2. $T' \cap \s$ consists of non-empty
collection of \sccs\ which are \ess\ in both $T'$ and $\s$.

\begin{proof}
Since $P'$ contains a punctured copy of $M$, $\s \cap P' \neq
\emptyset$. Note that the annulus $A$ is contained in the exterior
of $P'$. Since $\s$ is connected, $T' \cap \s \ne \emptyset$. By
the construction, it is clear that each component of $T' \cap \s$
is \ess\ in $T'$
(and, moreover, the slope of $T' \cap \Sigma$ is longitudinal
in $P' \cong (D^2 \times S^1) \# M$.)
Since the intersection $\s \cap T'$ can be
regarded as a subset of $\s \cap T$, we see that each component of
$T' \cap \s$ is \ess\ in $\s$.
\end{proof}

Let $X' = \text{cl}(M-P')$.

\medskip
\noindent
Subclaim 3.3.
The submanifold $X'$ satisfis the following.

\begin{enumerate}

\item
$X'$ is homeomorphic to the exterior of a non-trivial knot in $S^3$,
and

\item
there is a compressing disk, say $D_{X'}$, of $\partial X'$
such that $\partial D_{X'}$ is a meridian curve of $X'$.

\end{enumerate}

\begin{proof}
Recall that $T'$ bounds $P' \cong (D^2 \times S^1) \# M$. By the
construction, we see that $T'$ bounds a manifold homeomorphic to
$X \cup \text{cl}(\mx \setminus P)$ on the other side. Since the
slope of $\partial A$ is not meridional in $X$, we see that $X
\cup \text{cl}(\mx \setminus P)$ is homeomorphic to the exterior
of the same knot for $X$ (if the slope of $\partial A$ is
longitudinal) or, the exterior of a cable knot of the knot for $X$
(if the slope of $\partial A$ is not longitudinal). In either
case, the knot for $X'$ is non-trivial. It is clear that a
meridian disk for the solid torus factor of $P' \cong (D^2
\times S^1) \# M$ can be taken as $D_{X'}$.
\end{proof}

By Subclaims~3.2 and 3.3, we see that we may take $X'$ for $X$.
The procedures described above and in the proof of Claims~2 may repeated, if
necessary, and the process terminates since each application
reduces $| \s \cap T |$.
\end{proof}

\medskip
\noindent \underline{Claim 4.} $\smx$ compresses in $\mx$ into
both $H_1 \cup \mx$ and $H_2 \cup \mx$, and $\sx$ is essential in
$X$.

\begin{proof}
We first show the following:

\medskip
\noindent
Subclaim 4.1.
For each $i=1,2$,
there is a meridian disk $D_i$ of $H_i$ such that $D_i \cap T = \emptyset$.

\begin{proof}
Let $D \subset H_i$
be a meridian disk which minimizes $| D \cap T |$ among all
meridian disks. If $D \cap T = \emptyset$, then we are done.
Suppose that $D \cap T \ne \emptyset$. By the minimality of $| D
\cap T |$ no component of $D \cap T$ is a \scc. Then each
outermost disk in $D$ gives a boundary compression of $\sx
(\subset X)$ or $\smx (\subset \mx )$. Hence by Claims~2 and 3 and
Lemma~\ref{boudary compressible implies compressible}, we see that
there is a compressing disk $D_i$ for $\sx (\subset X)$ or $\smx
(\subset \mx )$ such that $D_i$ is contained in $H_i$,
however, this contradicts the minimality of $D$.
\end{proof}

Let $D_1$, $D_2$ be as in Subclaim~4.1. If one of $D_1$, $D_2$ is
contained in $\mx$, and the other in $X$, then the pair $\{ D_1,
D_2 \}$ gives a weak reduction for $\s$, a contradiction. Hence
either $\smx$ compresses into both sides or $\sx$ does.

\medskip
\noindent
Subclaim 4.2.
$\smx$ compresses into both sides.

\begin{proof}
Recall that $B = N(X \cup D_X)$. Then we minimize $| \s \cap
\partial B |$ via isotopy rel~$X$. If $| \s \cap \partial B | =
0$, then $\s$ is pushed into $B$, which is absurd. By using an
innermost disk argument, essentiality of the intersection and
irreducibility of $M$, we can show that each component of $\s \cap
B$ is essential in $\s$. Let $D^* (\subset \partial B)$ be an
innermost disk. Then $D^*$ is a meridian disk of $H_j$ $(j=1$ or
$2)$ contained in $\mx$. This together with the strong
irreducibility of $\Sigma$ shows that both $D_1$, $D_2$ are
contained in $\mx$. Hence $\smx$ compresses into both sides as
desired.
\end{proof}

By comment before Subclaim~4.2, we see that $\sx$ is
incompressible in $X$. By Claim~2,
we see that $\sx$ is not boundary parallel in $X$.
Hence $\sx$ is essential in $X$.
\end{proof}

\noindent
\underline{Claim 5.} $\smx$ is connected.

\begin{proof}
Since $\s$ is strongly irreducible, the compressions for $\smx$
(Claim~4) occurs on the same component of $\smx$. Assume that
there exists another component of $\smx$, say $F$. If $F$
compresses in $\mx$, then by the No Nesting Lemma (\ref{no
nesting}) the curve on $F$ that is compressed bounds a meridian
disk, say $D'$ in $H_j$ $(j=1$ or $2)$. Then $D'$ together with
$D_{3-j}$ gives a weak reduction for $\s$, a contradiction. Hence
$F$ is incompressible in $\mx \cong (D^2 \times S^1) \# M$. By 2
of Remark~\ref{remark of handlebody}, we see that $F$ is an
annulus, contradicting Claim~3. Hence $\smx$ is connected.
\end{proof}

Let $S_i'$ $(i=1,2)$ be a surface properly embedded in $\mx$ obtained
from $\smx$ by compressing into $H_i$-side maximally.

\medskip
\noindent
\underline{Claim 6.}
Each $S_i'$ is incompressible in $\mx$.

\begin{proof}
Assume that there is a compressing disk $D$ for $S_i'$ in $\mx$.
By isotopy of $D$ near its boundary we may assume that $\del D$ is
contained in $\smx$. By No Nesting Lemma(2.5) we may suppose that
$D$ is contained in $H_1$ or $H_2$. Since $S_i'$ is obtained from
$\smx$ by compressing into $H_i$-side maximally, $D$ must be
contained in $H_{3-i}$. However this shows that $\s$ is weakly
reducible, a contradiction.
\end{proof}

By Claim~6, 2 of Remarks~\ref{remark of handlebody}, and the
assumption of Step~1, each component of $S_1'$ (and $S_2'$) is an
annulus.

\medskip
\noindent \underline{Claim 7.} All the annuli of $T \cap H_i$
$(i=1,2)$ are simultaneously boundary compressible into $H_i \cap
\mx$, but not into $X$.

\begin{proof}
We will show the existence of boundary compressions into $\mx$.
This together with Claim 4, strong irreducibility and Lemma
\ref{boudary compressible implies compressible} then implies there
are no boundary compressions into $X$.

Since the argument is symmetric, it enough to prove Claim~7 for $T
\cap H_1$. Recall that $S_1'$ is an incompressible surface in $\mx
\cong (D^2 \times S^1) \# M$ obtained from $\smx$ by
simultaneously compressing into $H_1$ side. Hence the tubings for
retrieving $\smx$ from $S_1'$ are all done to the same side of
$S_1'$ and the tubes are not nested. Connectedness of $\smx$
implies that there is a unique $Z$, which is the closure of a
component of $\mx \setminus S_1'$, and within which all the
tubings are performed.  ($Z$ was obtained by attaching 2-handles
to $\mx \cap H_2$.)

By using innermost disk argument, we may assume that $S_1'$ and
the 2-sphere giving the connected sum structure $(D^2 \times S^1)
\# M$ are disjoint. For the analysis of the situation, we
temporarily ignore $M$ in $\mx \cong (D^2 \times S^1) \# M$. Then
each component of $S_1'$ is a boundary parallel annulus in $D^2
\times S^1$. We possibly have the following cases.

\medskip
\noindent
Case 1.
The components of $S_1'$ are not nested in $D^2 \times S^1$.

\medskip
\noindent
Case 2.
The components of $S_1'$ are nested in $D^2 \times S^1$.

\medskip
Suppose first that Case~1 holds.
If $S_1'$ is a single longitudinal annulus,
it boundary compresses into both sides, and
since retrieving $\smx$ is done via tubing into one side only,
Claim 7 clearly holds.  Else, let $Q$ be the union of the
parallelisms between the components of $S_1'$ and mutually
disjoint annuli in $T$. Then let $R = \text{cl}((D^2 \times S^1)
\setminus Q)$. Now recall punctured $M$ in $\mx$. We have the
following two subcases.

\medskip
\noindent Case 1.1. Punctured $M$ is contained in $R$ (\ie\ $Z =
R$).

\medskip
Note that $\smx$ is retrieved from the simultaneously boundary
parallel annuli $S_1'$ by adding tubes along
mutually disjoint arcs properly embedded in $R$.
This gives the conclusion of Claim~7.

\medskip
\noindent Case 1.2. Punctured $M$ is contained in $Q$ (\ie\ $Z
\subset Q$).

\medskip
Note that $\smx$ is retrieved from $S_1'$ by adding tubes along
mutually disjoint arcs properly embedded in $Q$. Hence by Claim~3,
we see that $S_1'$ consists of exactly one annulus (and $Z = Q$).
This implies that $R \cap T$ is an annulus properly embedded in
$H_1$. Denote this annulus by $A^*$.  By Claims~1, and 2 we see
that $A^*$ is an essential annulus in $H_1$. By
Lemma~\ref{essential annulus preserves spine}, there is a cycle of
a spine of $H_1$ in each side of $A^*$, in particular in $H_1 \cap
X$. This cycle is contained in $X$ and hence in the ball $B$, so
Frohman's Lemma (\ref{Frohman}) shows that $\s$ is reducible, a
contradiction. This shows that Case~1.2 does not occur.

\medskip
Suppose now that Case~2 holds (see Figure \ref{fig:case2}).

\begin{figure}[ht]

\begin{center}
\includegraphics[width=4cm, clip]{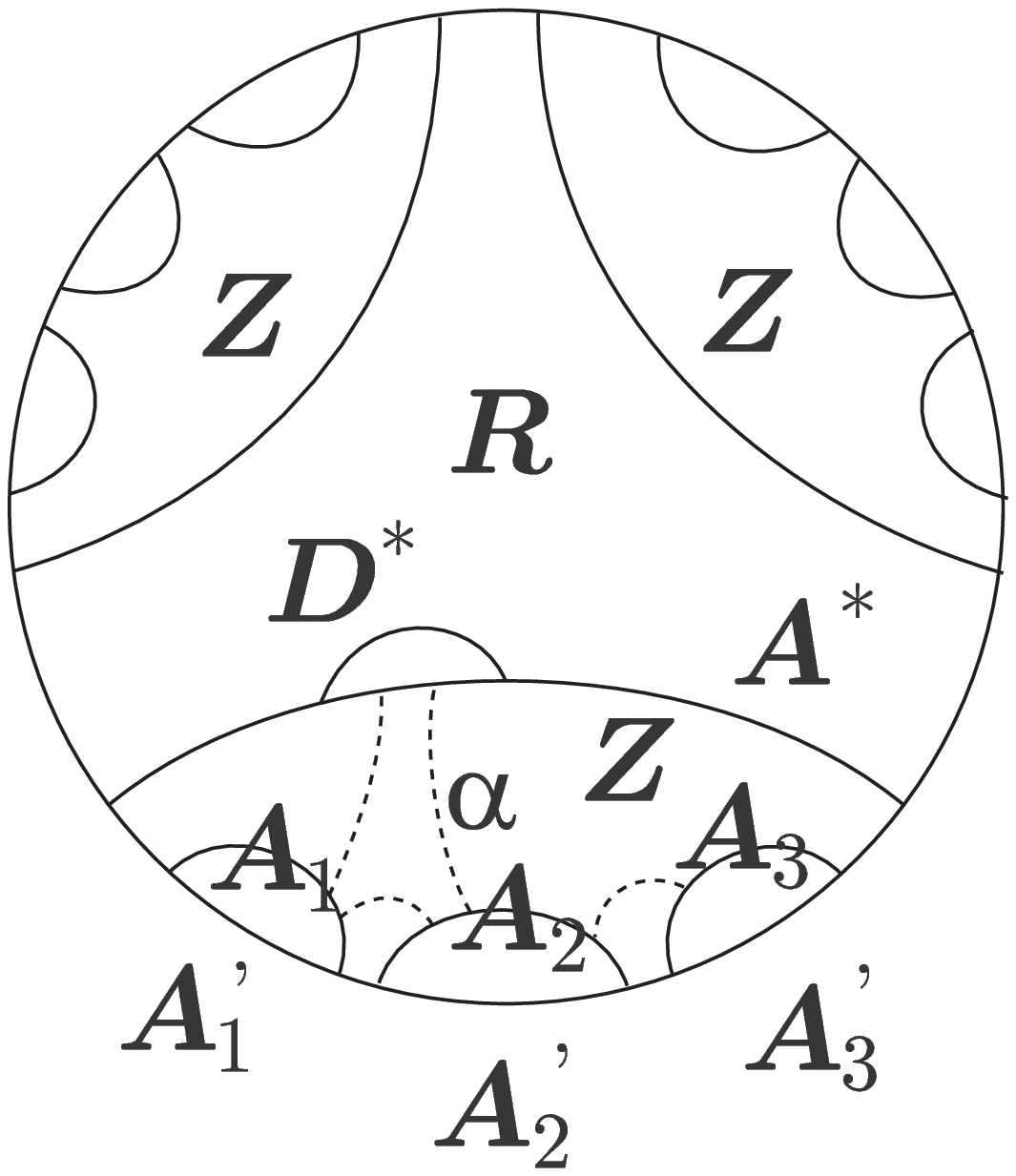}
\end{center}

\caption{Case 2.}
\label{fig:case2}
\end{figure}

Since $\sx$ is connected, and tubings for $S_1'$ for retrieving
$\sx$ is performed in one side of $S_1'$, we see that the depth of
the nesting is two (\ie\ there exist outermost and second
outermost annuli in $S_1'$, but no third outermost annulus), and
the tubings are performed along a system of mutually disjoint
arcs, say $\alpha$, properly embedded in the region between
outermost and second outermost components. Moreover connectedness
of $\smx$ implies that there exists exactly one second outermost
component, say $A^*$. Then the closures of the components of $(D^2
\times S^1)\setminus A^*$ consists of two components, say $P$ and
$R$, such that $P$ is a parallelism between $A^*$ and an annulus
in $T$, and $R$  a solid torus which contains a core of $D^2
\times S^1$. (Note that $Z \subset P$.) Again by the connectedness
of $\smx$, we see that every outermost component of $S_1'$, say
$A_1, \dots , A_n$, is contained in $P$.

Let $A_1', \dots , A_n'$ be annuli in $T$ such that
$\partial A_k'  = \partial A_k$, and
$A_k'$ and $A_k$ are parallel in $P$ $(k=1, \dots , n)$.
(Hence $A_1', \dots , A_n'$ are mutually disjoint.)
By simultaneously boundary compressing $A_1', \dots , A_n'$
into $\smx$ in $P$, we obtain disks, say $D_1, \dots , D_n$,
properly embedded in $H_1$. Let $D^*$ be a disk properly embedded
in $R$ such that $\partial D^*$ is a \scc\ in $A^*$ which bounds a
disk (in $A^*$) containing the points $\alpha \cap A^*$. We may
regard $D^*$ as a disk properly embedded in $H_1$, and it is clear
from the definition that $D_1 \cup \dots \cup D_n \cup D^*$ cuts
off a handlebody corresponding to $H_1 \cap \{$the region between
$A_1 \cup \cdots \cup A_n$ and $A^*\}$. Note that the exterior of
this handlebody in $H_1$ is a non-trivial handlebody, since it contains
\sccs\ $\partial A^*$ which are essential in $\s$. This shows that
there is a cycle of a spine of $H_1$ that is contained in the
exterior of the punctured $M$. Hence by Frohman's Lemma
(\ref{Frohman}) we see that $\s$ is reducible, a contradiction.
This shows that Case~2 does not occur.

This completes the proof of Claim~7.
\end{proof}

\noindent
\underline{Completion of the proof of Step~1.}

\medskip
By Claim~7, we see that $T \cap H_1$ consists of annuli that are
simultaneously boundary compressible into $\mx$ side. Let $K$ be the
closure of a component of $H_1 \setminus T$ to which the boundary
compressions are \em not \em performed. By Claim~7, we see that
$K$ is contained in $X$. By performing the boundary compressions
on $T \cap H_1$, we obtain a union of mutually disjoint meridian
disks (say $\widehat{D}$) in $H_1$. Let $K'$ be the closure of the
component of $H_1 \setminus \widehat{D}$ such that $K' \supset K$.
Since we obtain $A^*$ from $\widehat D$ by banding into $K'$, we
see that $K'$ is not a ball.
Hence $K'$
contains a cycle of a spine of $H_1$, and $K$ contains the same
cycle. By Frohman's Lemma (\ref{Frohman}), we see that $\s$ is
reducible. This contradiction completes the proof of Step~1.

\medskip
\noindent
{\bf Step 2: $T \cap H_j$ $(j=1$ or $2)$ contains a boundary parallel annulus.}

\medskip
Recall that the assertion of Step~2 is conclusion 1 of
Theorem~\ref{main theorem}.

\medskip
\noindent \underline{Claim 8.} There is an annulus component of $T
\cap H_j$ $(j=1$ or $2)$ which is compressible in $H_j$.

\begin{proof}
Let $\gamma$ be a component of $T \cap \s$. Since $\gamma$ defines
a meridional slope (Step~1), it bounds a disk such that a \nbhd\
of $\gamma$ in the disk is embedded in one of the handlebodies
$H_1$ or $H_2$. By the No Nesting Lemma (\ref{no nesting}),
$\gamma$ bounds a disk that is entirely in $H_1$ or $H_2$. By
innermost disk argument applied to the intersection of this disk
with $T$, we see that some annulus of $T\setminus \s$ compresses
in some $H_j$.
\end{proof}

Let $A$ be the annulus obtained in Claim~8. Without loss of
generality, we may suppose that $A \subset H_1$. Let $A'$ be an
annulus component of $T \cap H_2$ adjacent to $A$, (\ie , $A$ and
$A'$ share a common boundary component, say $\gamma$). Since $A$
is compressible in $H_1$ (Claim~8), $\gamma$ bounds a disk in
$H_1$. If $A'$ compressed in $H_2$, then $\gamma$ bounds a disk in
$H_2$ too, and this shows that $H_1 \cup_\s H_2$ is reducible, a
contradiction. Hence $A'$ is incompressible in $H_2$. Assume that
$A'$ is not boundary parallel in $H_2$. Since $A'$ is
incompressible, $A'$ is boundary compressible (1 of
Remark~\ref{remark of handlebody}). Let $D'$ be a disk obtained
from $A'$ by a boundary compression. Since $A'$ is incompressible
and not boundary parallel, $D'$ is a meridian disk. By applying a
slight isotopy, we may suppose that $\partial D' \cap \gamma =
\emptyset$. Hence $D'$ together with a disk obtained by
compressing $A$ gives a weak reduction of $H_1 \cup_\s H_2$, a
contradiction.

This completes the proof of Step~2.

\medskip
\noindent
{\bf Step 3: Completion of the proof.}

\medskip
Finally we induct on $|T \cap \s|/2$ to show that every component
of $\sx$ is a meridional annulus. Let $A' (\subset H_j)$ be an
annulus obtained in Step~2, and $A_\s$ the annulus in $\s$ such
that $A \cup A_\s$ bounds a parallelism $P$ contained in $H_j$.

Suppose that $|T \cap \s|/2 = 1$. Then either $\sx = A_\s$ or
$\smx = A_\s$. However if $\smx = A_\s$, then $\s$ can be isotoped
into the 3-ball $B$, a contradiction. Hence $\sx = A_\s$, which
gives the conclusion~2 of Theorem~\ref{main theorem}.

Suppose that $|T \cap \s|/2 > 1$. By passing to outermost one, if
necessary, we may suppose that $\text{int} P \cap \s = \emptyset$.
We have the following two cases.

\medskip
\noindent
Case 1.
$A_\s \subset X$.

\medskip
In this case, we push $A_\s$ along the parallelism $P$ out of $X$.
Then by induction, we see that the image of $\s$ intersects $X$ in
meridional annuli. Note that $\sx$ is the union of these annuli
and $A_\s$. Hence each component of $\sx$ is an annulus, and their
slope is meridional by Step 1.

\medskip
\noindent
Case 2.
$A_\s \subset \mx$.

\medskip
In this case, we push $A_\s$ along the parallelism $P$ into $X$.
Then by induction, we see that the image of $\s$ intersects $X$ in
meridional annuli. To retrieve $\sx$ we push a core curve of one
of these annuli out of $X$. Thus this annulus breaks into two
annuli. Again each component of $\sx$ is a meridional annulus.

This completes the proof of Theorem~\ref{main theorem}.
\end{proof}

\end{document}